# Risk Assessment for Performance Driven Building Design with BIM-Based Parametric Methods

Fatemeh Shahsavari, Jeffrey D. Hart, Wei Yan
Texas A&M University, College Station TX 77843, USA
Fatemeh.shahsavari@tamu.edu
hart@stat.tamu.edu
wyan@tamu.edu

**Abstract.** A growing demand for handling uncertainties and risks in performance-driven building design decision-making has challenged conventional design methods. Thus, researchers in this field lean towards viable alternatives to using deterministic design methods, e.g., probabilistic methods. This research addresses the challenges associated with conventional methods of performance-driven building design, i.e., ignoring the existing uncertainties and lacking a systematic framework to incorporate risk assessment in building performance analysis. This research introduces a framework (BIMProbE) to integrate BIM-based parametric tools with building probabilistic performance analysis to facilitate uncertainty analysis and risk assessment in performance-based building design decision-making. A hypothetical building design scenario is used to demonstrate the application of the proposed framework. The results show that the probabilistic method leads to a different performance ranking order than the deterministic method. Also, the probabilistic method allows evaluating design options based on different risk attitudes.

## 1 Introduction

Architectural design decision-making begins with identifying design problems and objectives. It sets boundaries for designers' potential problem-solving methods. Performance-based building design focuses on methods and strategies that integrate and optimize different aspects of building performance. Using the computer power allows designers to explore a broader range of solutions, efficiently. Building performance simulation (BPS) tools have been extensively used by the architects and engineers to simulate building performance [1].

The BPS tools usually produce building performance predictions based on a set of input data including building physical characteristics, interior conditions, weather data, and mechanical specifications [2], [3], which commonly come with outstanding uncertainties [4]. The Building Energy Software Tools Directory provides an ongoing list of 200 simulation tools [5], but most of these tools are not designed to deal with uncertainties. The majority of the BPS tools only collect deterministic input data and run a single deterministic simulation to evaluate the building performance [6], [7].

The various sources of uncertainties associated with building performance and its potential impact on the human health and environmental crisis demands for new methods of design thinking to deal with uncertainties. The advancement of computer technologies allows the application of data-driven methodologies in the field of building performance analysis to achieve resiliency and sustainability.

This research introduces a framework to incorporate probabilistic models into the building design decision making process, demonstrated with a design test case. The goal of this research is to tackle data uncertainties and potential risks in architectural design decision-making with a focus on building energy performance. For that purpose, the results obtained from the proposed probabilistic framework are compared with those from a conventional deterministic method. Also, three design decision making criteria including expected value, maximax, and maximin are applied to discuss the simulation results based on different attitudes towards risk.

# 2 Literature Review

Deterministic methods in performance-driven building design decision-making fail to address existing uncertainties in design decision making. Uncertainty analysis techniques such as Monte Carlo coupled with risk assessment methods introduce a potential solution to tackle uncertainties and make robust design decisions [8].

## 2.1 Monte Carlo and building performance simulation

As Saltelli et al. [9] stated, variance-based methods, e.g., Monte Carlo, have shown more effectiveness and reliability when working with uncertainties. The Monte Carlo simulations use random variables and input probability density functions to address the stochastic status of the problem. Let a mathematical modeling $Y = f(x)$ define correlations between a vector of input variables $X = \{X_1, X_2, \ldots, X_k\}$ and an output $Y$ , where $f$ is a deterministic integrable function which translates from a $k - D$ space into a $1 - D$ one, i.e., $\mathbb{R}^K \rightarrow \mathbb{R}$. The model produces a single scalar output $Y$ when all input variables are deterministic scalars. However, if some inputs are uncertain or undecided, the output $Y$ will also associate with some uncertainties. An input variable $X_i,$ is defined by a mean value $\mu_i$, a variance $\sigma_i$, and a probability distribution, such as Normal, Uniform, Poisson, etc. In the Monte Carlo methods, a set of samples from possible values of each input variable are generated. These input values are inserted into the simulation model to generate the probability distribution of the output $Y$. Processing the output range $Y$ delivers a mean value and a frequency distribution of the output.

The challenge of mapping between simulation tools and probabilistic techniques in the process of Monte Carlo simulations has received a lot of attention in the literature. Lee et al. [10] introduced an uncertainty analysis toolkit explicitly for building performance analysis referred to as the Georgia Tech Uncertainty and Risk Analysis Workbench (GURA_W). The identification and modification of input variables are possible using GURA_W, and the uncertainty quantification repository available in this tool allows the energy modelers to access the uncertainty distributions of previous parameters being modeled. Hopfe et al. [6] studied uncertainty analysis in an office building energy consumption and thermal comfort assessment through connecting MATLAB with a building performance simulation tool, called VA114. They declared that including uncertainties in building performance analysis could support the process of building design decision making. Macdonald and Strachan [11] integrated uncertainty analysis with building energy simulation using Esp-r software. They studied the uncertainties in thermal properties of building materials and building operation schedules and concluded that uncertainty analysis facilitates risk assessment and improves building design decision making. de Wit and Augenbroe [12] studied the impact of variations in building material thermal properties, along with model simplifications, on building thermal behavior, using two simulation tools, ESP-r and BFEP, integrated with Monte Carlo uncertainty analysis technique. They applied expected value and expected utility decision making criteria to discuss how designers can use the extra information obtained from uncertainty analysis. Asadi et al. [13] used Python programming to automate the integration of Monte Carlo simulations into energy analysis. They developed a regression model as a pre-diagnostic tool for energy performance assessment of office buildings to identify the influence of each design input variable.

### 2.2 BIM and parametric tools integrated with uncertainty analysis and risk assessment

BIM technologies may be useful in performance analysis and risk management since they facilitate transferring data from BIM authoring tools to other analysis tools, also enable designers to automate iterations in design and analysis processes [14], [15]. Kim et al. [16] applied the Monte Carlo uncertainty analysis in building energy analysis through the integration of BIM and MATLAB platforms. A set of software applications including Revit Architecture 2010, ECOTECT 2010, and EnergyPlus 6.0 were used for modeling and simulation. For uncertainty analysis, the MATLAB Graphical User Interface (GUI) platform was applied to develop a self-activating Monte Carlo simulation program. Rezaee et al. [17] provided a CAD-based inverse uncertainty analysis tool to estimate the unknown input variables and improve design decision-makers confidence in the early stage of building design. They created two energy models, one in EnergyPlus and another one in a spreadsheet-based energy analysis tool to run the energy calculations.

Although BIM and parametric design tools allow the iteration of building performance simulations, the integration of Monte Carlo into the platform of these applications have not been broadly studied. For the performance-based building design, further studies are required to provide clear guidance on the mapping between BIM authoring tools, parametric analysis tools, and probabilistic techniques such as Monte Carlo to provide probabilistic outcomes for building energy analysis [14]. This research intends to apply variance-based methods such as Monte Carlo in the performance driven design decision making process using BIM and parametric tools.

## 2.3 Design decision making and risk assessment

The international ISO standard (ISO 31000:2009) defines "risk" as the effect of "uncertainty" on objectives, that could be positive consequences as well as negative impacts [18]. Different decision-making criteria including expected value, maximax, and maximin for risk assessment and design decision making under uncertainties are applied in this research. In this approach, two major terms are introduced:

- Key performance indicators (KPI) [19], and
- Key risk indicators (KRI) [20].

While KPI is a metric to measure the performance or objective, KRI measures deviations from the target and depicts the threat of breaching a specific threshold [21]. KPI and KRI are useful metrics to measure the objective aspect of risk (likelihood), also facilitate decision making based on the risk attitude of decision makers and evaluate the performance subjectively.

A building energy model is too complex to estimate the variance of the output by just reviewing the variance of the input variables. In building energy projects, KPI could be defined as building annual energy consumption or energy saving, and KRI may show the possibility of occurrence of discrepancy between expected energy performance during the design stage and real energy performance after project completion [22]. Thus, performing a probabilistic analysis and uncertainty propagation for identifying the probability distribution of the output could be helpful in the process of risk assessment and design decision making.

**Expected value criterion**

The Expected Value (EV) of a variable is the return expected or the average benefit gained from that variable. This statistical measure, the sum of all possible gains, each multiplied by their probability of occurrence, demonstrates the cost-benefit analysis of a design option, considering the input uncertainties [23].

De Wit and Augenbroe [12] applied the criterion of EV to evaluate competing design alternatives. Their example describes a situation where a designer needs to make a decision on whether or not to use a cooling system in their building design. The designer would decide to use a cooling system only if the indoor temperature excess of the building without a cooling system is more than 150 hours. Based on their probabilistic results, they concluded that the most likely value of the outcome or EV would be well below 150 hours and the designer could make their decision comfortably not to include the cooling system in their design.

**Maximin and maximax criteria**

Design options are not always assessed objectively, as suggested by the EV theory. It is important to know how likely different design outcomes are to occur for a design option, but decision-makers' preferences affect their decision in selecting the optimal design option as well [23]. The best choice for one design decision-maker might not be the best for another one with different preferences, thus it is not always the best decision to select the design option with the maximum EV. There may be a situation which demands risk-averse decision-making and selecting a design option which can go the least wrong. This decision-making strategy is known as maximin and suggests selecting the design option maximizes the minimum payoff achievable [24], [25]. On the other hand, taking aa risk-seeking approach towards risk might lead to selecting a design option with the most optimistic possible outcome. This decision-making strategy is termed maximax and searched for a design option to maximize the maximum payoff available [24], [26].

## 3 Methodology

This study presents a new framework to implement probabilistic methods in the field of performance-driven building design decision-making, using Building Information Modeling (BIM) and parametric tools. Fig. 1 presents the workflow including the key steps of this process, data being transferred, variables, and the software:

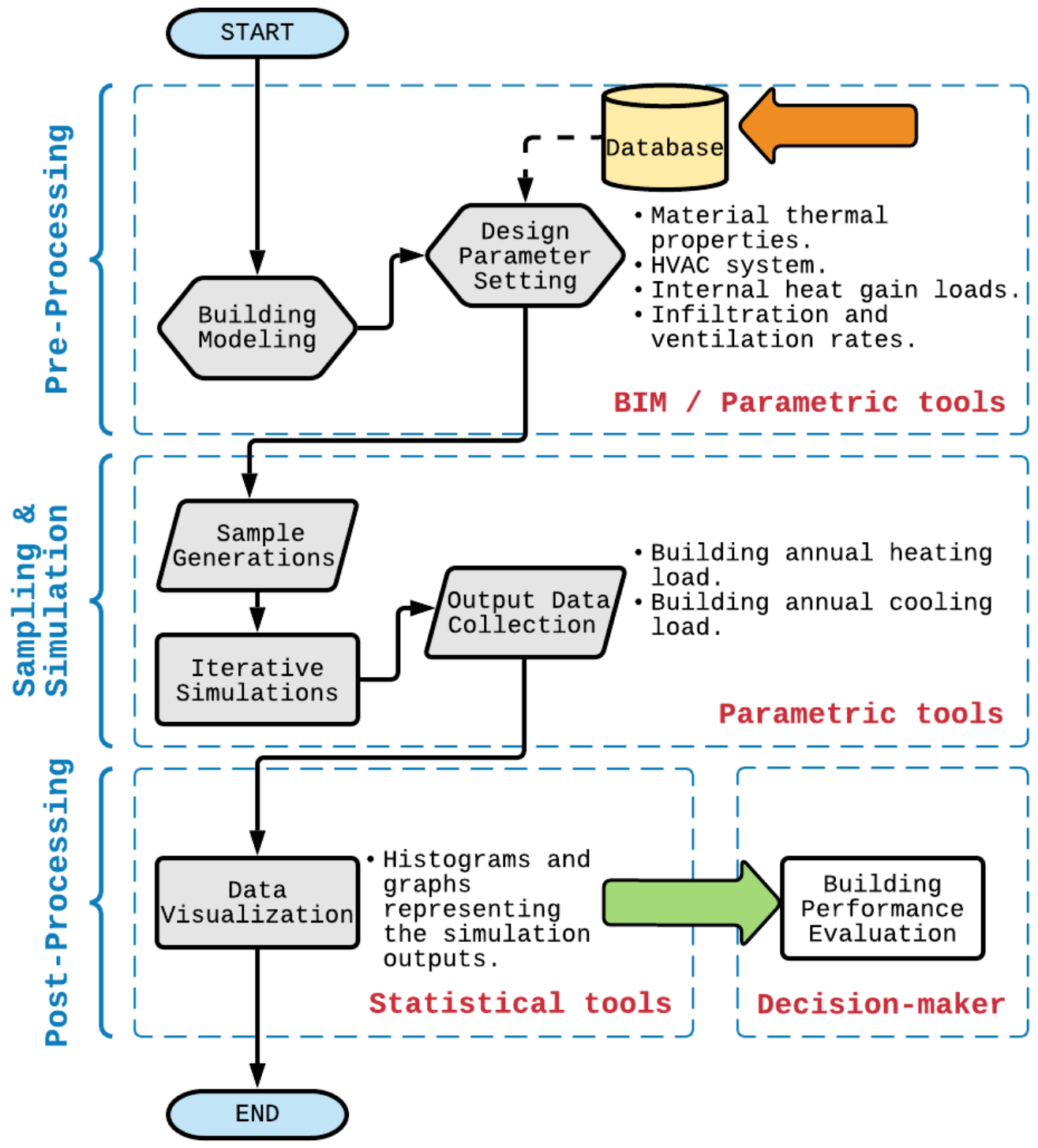

**Fig. 1.** BIMProbE workflow for probabilistic performance-based design decision-making.

This framework considers the uncertainties in building energy simulations including material properties, internal heat gains, and infiltration and ventilation rates (more details about this framework are presented in [27]). This framework targets the Revit-based building design process and has three main steps described below:

### 3.1 Pre-processing

The building geometry is modeled in a BIM authoring tool such as Revit. An external application named BIMProbE is developed to set the probability distributions for thermal properties of building materials based on the real-world material properties learned from literature studies and retrieved from an external database (Microsoft Excel). The mapping between the BIM model and the Excel-based database is performed through BIM API and using object-oriented programming.

Building elements consist of exterior walls, interior walls, roofs, floors, and windows. This program will find and collect opaque components including walls, roofs, and floors. The material IDs associated with each material is identified to collect the thermal properties. The thermal properties of building materials including thickness, thermal conductivity, and R-value can be collected from the BIM model. For instance, the thickness can be accessed as the width of each layer of the building components. As an example, a Structurally Insulated Panel (SIP) wall type consists of six layers of Plasterboard, Wood, Timber Insulated Panel – OSB, Timber Insulated Panel – Insulation, Timber Insulated Panel – OSB, and Sand/Cement Screed. The width of each layer of an SIP wall can be identified as the thickness of that specific material. If the thermal properties of a material are missing in the Revit model, the program will create a thermal asset for that material and will set the thermal properties according to the corresponding values in the Excel database.

BIMProbE add-in will automatically create the required shared parameters to add the probability distributions to the building material thermal properties and bind them to each building material in the project. This application searches for the material names in the Excel database to assign the correct value to each parameter. This process begins with an attempt to open the external Excel database. Once the user starts this addin, a window will pop up asking the user to select an external Excel file. If the user selects an Excel file, the program will automatically start reading it cell by cell. The program will check the sanity of the data first and if there is no error found, it will continue with the rest of this process. The probability distributions (including the mean and standard deviation) for each material type are set according to the corresponding values in the Excel database. At the end of this process, BIMProbE allows writing the mean and standard deviation values of the thermal properties to a new Excel spreadsheet to be used for later steps of sampling and simulation.

Other design input parameters required for building energy simulations including Heating, Ventilation, and Air Conditioning (HVAC) system specifications, internal heat gain loads, infiltration and ventilation rates, and operational schedules are added using parametric tools such as Grasshopper, which is the visual programming environment for Rhinoceros.

### 3.2 Sampling and simulation

The information prepared in the pre-processing step is used for sample generation and iterative simulations. Using Latin Hypercube Sampling method, N samples for each design input variable are generated and N energy simulations are run for each design option (building forms) in Grasshopper. The simulation outputs (building annual heating and cooling loads) are recorded in Grasshopper and written to Excel using TT Toolbox plugin for Grasshopper [28]. For each round of the energy simulation, a corresponding value from the sample pool of each design variable ($x_{uncer_i}$) is selected and inserted into the simulation model, and as a result, N output values are obtained. The set of outputs whose elements correspond to the samples are $Y = \{y_1, y_2, y_3, \dots, y_m\}, m = N$, where $Y$ denotes the space of output values, which are the results of building thermal load simulations.

In the deterministic approach, all the design variables, regardless of being deterministic or uncertain, are assigned to their associated mean ($\mu_i$) as a fixed value ($x_i = \mu_i$). Fixing all the input variables at their mean values, the energy simulation is run once and a single output $y_{det} = f(\mu_1, \mu_2, \mu_3, \dots, \mu_n)$ is obtained, where $(\mu_1, \mu_2, \mu_3, \dots, \mu_n)$ are the means of the design variables.

### 3.3 Post-processing and design decision-making

The post-processing phase consists of data analysis and graphical presentation of the simulation results. This phase is conducted using a statistical software known as JMP. The data are collected in Grasshopper and exported to Excel for post-processing. The JMP [29] add-in for Excel provides interactive graphics and tables that enable the user to identify relationships visually and examine patterns. The histogram demonstrations, normality plots and box plots are used for risk assessment and data visualization.

The output results are represented by the values of deterministic and probabilistic outputs using two metrics of KPI and KRIs. In addition, design options are ranked according to the criteria of expected value, maximax, and maximin. The information provided in this process could help design decision-makers with building performance evaluation including uncertainties and comparing the probabilistic results with the deterministic results.

## 4 Test Case

This test case studies the energy performance of a hypothetical mid-size office building with four design options. For this test case, the building annual thermal load (the sum of heating and cooling loads) in kWh/m$^2$ is defined as the target. The effects of some of the existing uncertainties including the thermal properties of building materials, internal heat gain loads, ventilation rates, infiltration rate, and occupant behavior are studied.

### 4.1 Test case model

This test case presents a hypothetical 5-story office building, 54m by 93m (Total gross area of 25,110 m$^2$) in Chicago, Illinois (climate zone 5A). Fig. 2. shows the 3D view of the building base model in Revit, and the building visualization in Rhino using Ladybug tools in Grasshopper, from left to right. The building geometry is created in Revit and the materials are defined in the same model. The building information including the materials' thermal properties is exported from Revit to Excel. This information is used in Rhino/Grasshopper to generate samples and run multiple energy simulations using the parametric capabilities of Grasshopper.

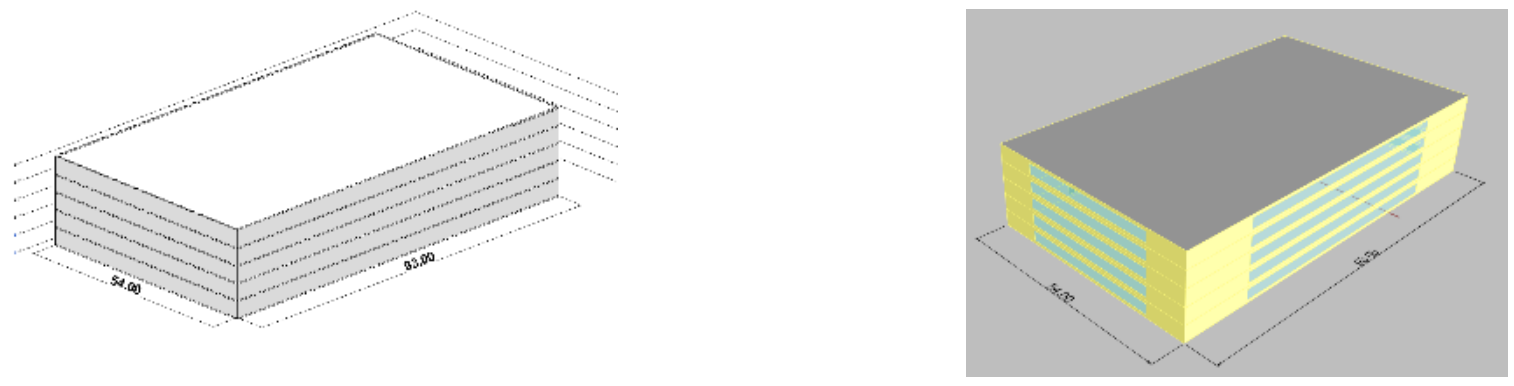

**Fig. 2.** Geometric representation of the test case.

The building geometry, layout, orientation, operation schedules, internal heat gain sources, and HVAC type are assumed to be identical throughout this experiment. On the other hand, exterior wall construction, floor construction, roof construction, Window-to-Wall Ratio (WWR), glazing properties, and temperature set points vary among the proposed design options (see Table 1).

**Table 1.** Design variable assumptions.

| Input Variable | Design Option 1 | Design Option 2 | Design Option 3 | Design Option 4 |
|---|---|---|---|---|
| RSI of Exterior wall construction [$m^2K/W$] | 3.7 | 2.76 | 3.35 | 3.35 |
| RSI of Floor construction [$m^2K/W$] | 4.59 | 4.59 | 4.59 | 3.35 |
| RSI of Roof construction [$m^2K/W$] | 5.88 | 4.17 | 5.88 | 5.88 |
| Window to Wall Ratio (WWR)% | 70/65/60/20 (N/S/ E/W) | 40 overall | 50/70/40/45 (N/S/ E/W) | 50/70/40/45 (N/S/E/W) |
| Glazing U-value [$W/m^2K$] | 1.7 | 3.12 | 1.2 | 1.2 |
| Glazing SHGC | 0.20 | 0.42 | 0.20 | 0.15 |

The materials in this test case are specified with some attributes such as RSI value and U-value. These terms are commonly used in the building industry to describe the thermal resistance and thermal conductivity of building materials. The mean values of these variables are defined based on the ASHRAE 90.1-2010 requirements for this test case's climate zone.

The internal heat gain loads, ventilation rates, and infiltration rate, along with the HVAC system description and building operation schedules are used to convert the building mass into thermal zones. The realizations of the mean values for the internal heat gain loads, ventilation rates, and infiltration rate for all four design options are specified based on the ASHRAE 90.1-2010 requirements for climate zone 5A (Table 2. ).

**Table 2.** Description of the system design assumptions.

| **Variable [Unit]** | **μ** |
|---|---|
| Equipment loads per area [W/m$^2$] | 10.765 |
| Infiltration (air flow) rate per area [m$^3$/s-m$^2$] | 0.0003 |
| Lighting density per area [W/m$^2$] | 10.55 |
| Number of people per area [ppl/m$^2$] | 0.07 |
| Ventilation per area [m$^3$/s-m$^2$] | 0.0006 |
| Ventilation per person [m$^3$/s-m$^2$] | 0.005 |

The standard deviations of all uncertain variables are assumed to be equal to 10% of the corresponding mean values, due to a lack of information. This assumption is based on the previous research done by [30]. de Wit [30] has estimated this percentage to be up to 10% in his report. For the future research, it is encouraged to conduct experiments or use expert knowledge to provide an appropriate standard deviation for each input variable with uncertainties.

The other system design inputs are assumed to be fixed and described as shown in Table 3.

**Table 3.** HVAC system specifications

| | **Design Option 1** | **Design Option 2** | **Design Option 3** | **Design Option 4** |
|---|---|---|---|---|
| HVAC system type | VAV with reheat | VAV with reheat | VAV with reheat | VAV with reheat |
| Temperature set points for cooling/ heating | 23˚C/ 21˚C | 23˚C/ 21˚C | 23˚C/ 21˚C | 23˚C/ 21˚C |
| Supply air temperature for cooling/ heating | 12.78˚C/ 32.2˚C | 12.2˚C/ 32.2˚C | 12.78˚C/ 32.2˚C | 12.2˚C/ 32.2˚C |
| Chilled water temperature | 7.2˚C | 6.7˚C | 7.2˚C | 6.7˚C |
| Hot water temperature | 60˚C | 82.2˚C | 60˚C | 82.2˚C |
| maximum heating supply air temperature | 40˚C | 40˚C | 40˚C | 40˚C |
| minimum cooling supply air temperature | 14˚C | 14˚C | 14˚C | 14˚C |
| maximum heating supply air humidity | 0.008 kg-$H_2O$/kg-air | 0.008 kg-$H_2O$/kg-air | 0.008 kg-$H_2O$/kg-air | 0.008 kg-$H_2O$/kg-air |
| ratio | | | | |
| minimum cooling supply air humidity ratio | 0.0085 kg-$H_2O$/kg-air | 0.0085 kg-$H_2O$/kg-air | 0.0085 kg-$H_2O$/kg-air | 0.0085 kg-$H_2O$/kg-air |
| recirculated air per area | 0 $m^3$/s-$m^2$ | 0 $m^3$/s-$m^2$ | 0 $m^3$/s-$m^2$ | 0 $m^3$/s-$m^2$ |

The occupancy, lighting, and equipment schedules are matched with the office schedules in the ASHRAE 90.1-2010 [31] for all four design options.

## 4.2 BIMProbE and probabilistic energy analysis

The building element modeling environment in Revit is used to model the base model, assign the associated materials, and add the probability distributions to the thermal properties of building materials. BIMProbE add-in for Revit allows creating required shared parameters to add probability distributions of the thermal properties of building materials and export this information to Excel for further analysis. Fig. 3 demonstrates the Revit model and the workflow to set the probability distributions for thermal properties of building materials.

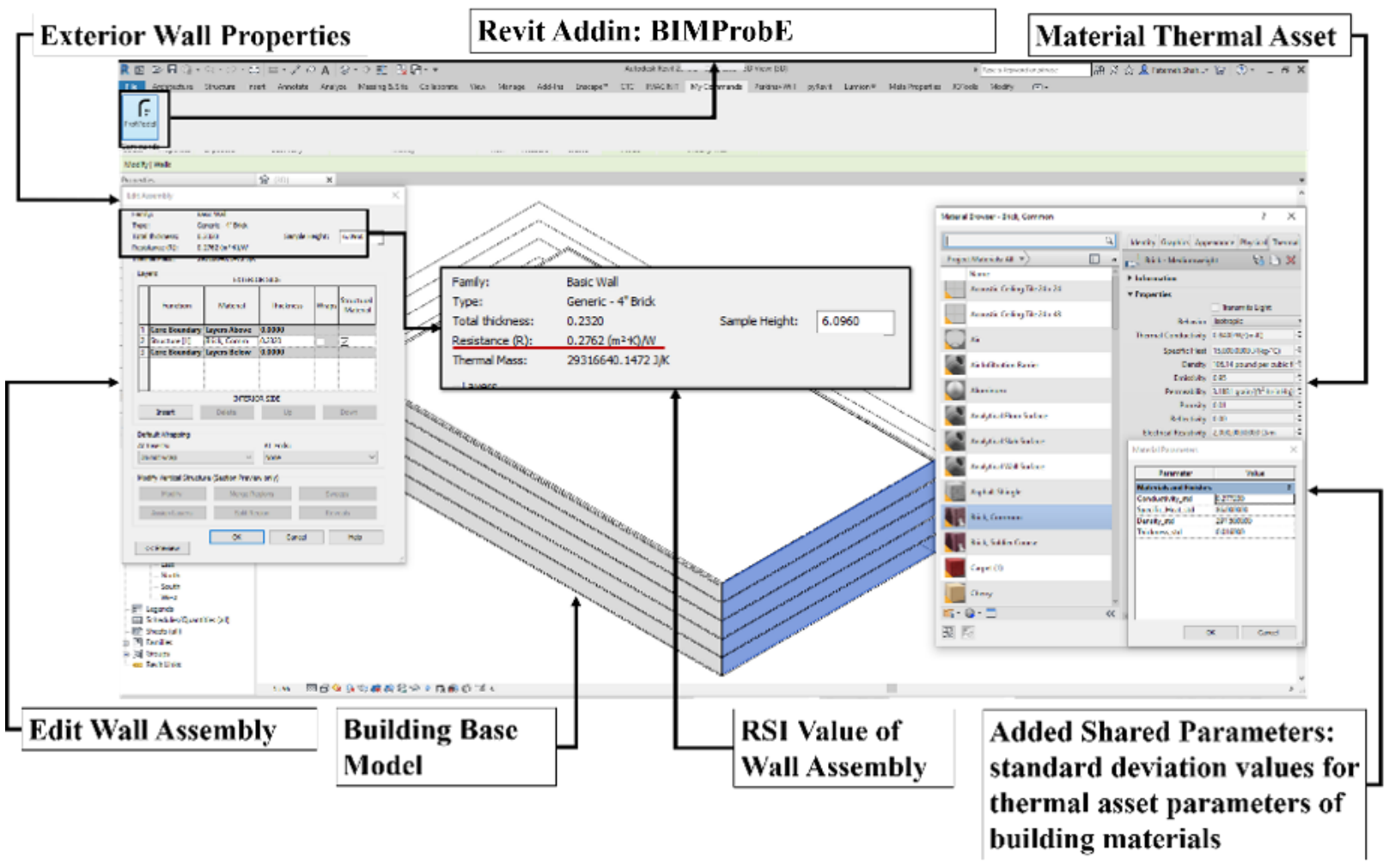

**Fig. 3.** BIMProbE add-in for getting and setting the probability distributions of thermal properties for building materials.

The visual programming tool in Rhino, Grasshopper, is used to get the base model from Revit and develop the four different design options. The building mass was converted to thermal zones in each case with defining the adjacency types, WWRs, internal heat gain loads, and HVAC system specifications. The building mass is split to floors and thermal zones using the existing components in Ladybug tools. The WWR for each façade is parametrically set in Grasshopper.

In this test case, the uncertain variables for which the variations take place due to unpredictable changes during construction, climate change, age, and maintenance are sampled using normal probability distribution. On the other hand, bases on the findings of [32] the input variables related to occupant behavior or presence can be best described with Poisson distributions.

The mean and standard deviation values of the thermal properties of building materials are set according to the Excel data inventory (created using BIMProbE).

The probabilistic simulations are conducted in Grasshopper using the statistical tools programed in CPython and simulation applications available in Ladybug tools (with OpenStudio simulation engine). The CPython component in Grasshopper is used to import the statistical tools such as Numpy and Scipy into Grasshopper (Abdel Rahman, 2018) to generate input samples with normal and Poisson distributions.

The internal heat gain sources and the HVAC settings are set in Grasshopper based on the findings of the previous research. The occupancy, lighting, and equipment operation schedules are set, and the thermal zones are exported to IDF files and run through OpenStudio in Grasshopper. The energy simulations are programmed to start automatically and run using the generated input samples. To automate the random value selection for each input variable and run OpenStudio for 500 (number of samples) times, a number slider, which is controlled by the Fly component in Ladybug tools, is connected to the list of input variables and selects an index of each list (starting from 0 and ending at 499) automatically and feeds the associated input value to the simulation. The simulation output is the building annual thermal load calculated in kWh/m$^2$. The outputs are collected and stored in Excel for post-processing.

The design variables with uncertainties are denoted as $X_{uncer} = \{x_{uncer_1}, x_{uncer_2}, x_{uncer_3}, \dots, x_{uncer_k}\}$, k = 11. In this test case, the eleven input variables with uncertainties include the RSI values of exterior walls ($x_{uncer_1}$), the glazing U-value ($x_{uncer_2}$), the RSI values of floor construction ($x_{uncer_3}$), the RSI values of roof constructions ($x_{uncer_4}$), the three internal heat gain loads (equipment ($x_{uncer_5}$), lighting ($x_{uncer_6}$), and people loads ($x_{uncer_7}$), the infiltration rate ($x_{uncer_8}$), the two ventilation rate factors (ventilation per person ($x_{uncer_9}$) and ventilation per area ($x_{uncer_{10}}$), and the infiltration schedule that is strongly correlated with the possibility of opening or closing windows by occupants ($x_{uncer_{11}}$) for which the mean and standard deviation values are obtained from the literature studies [33], [34]. In this test case, the probability distributions of physical input variables $\{x_{uncer_1}, x_{uncer_2}, x_{uncer_3}, x_{uncer_4}\}$, also ($x_{uncer_8}$), ($x_{uncer_9}$), and ($x_{uncer_{10}}$) are set as normal, while the internal heat gain loads $\{x_{uncer_5}, x_{uncer_6}, x_{uncer_7}\}$, also the infiltration schedule ($x_{uncer_{11}}$) are assumed to be highly dependent on the occupant presence and sampled using the Poisson distribution.

# 5 Results

In this research, the main design objective is improving building energy performance, thus KPI and KRI are defined to rank the predicted energy performance of different design options. KPI in this study is building annual thermal (heating and cooling) loads, and the lowest possible value would be desired. KRIs include the mean, standard deviation, and variance of building annual thermal loads. The mean value shows the average of the samples, while standard deviation and variance are presented as additional measures of risk. The probabilistic framework works with quantifying the uncertainties in design inputs and allows predicting the probability distribution of the simulation outcome and the risks threatening building energy performance for each design option.

Using a processor Intel® Core™ i7-4770 CPU at 3.40 GHz speed, 500 energy simulations were run for each design option, each simulation taking 11 seconds. Fig. 4. illustrates deterministic and probabilistic results of annual thermal loads (kWh/m$^2$) for the four design options.

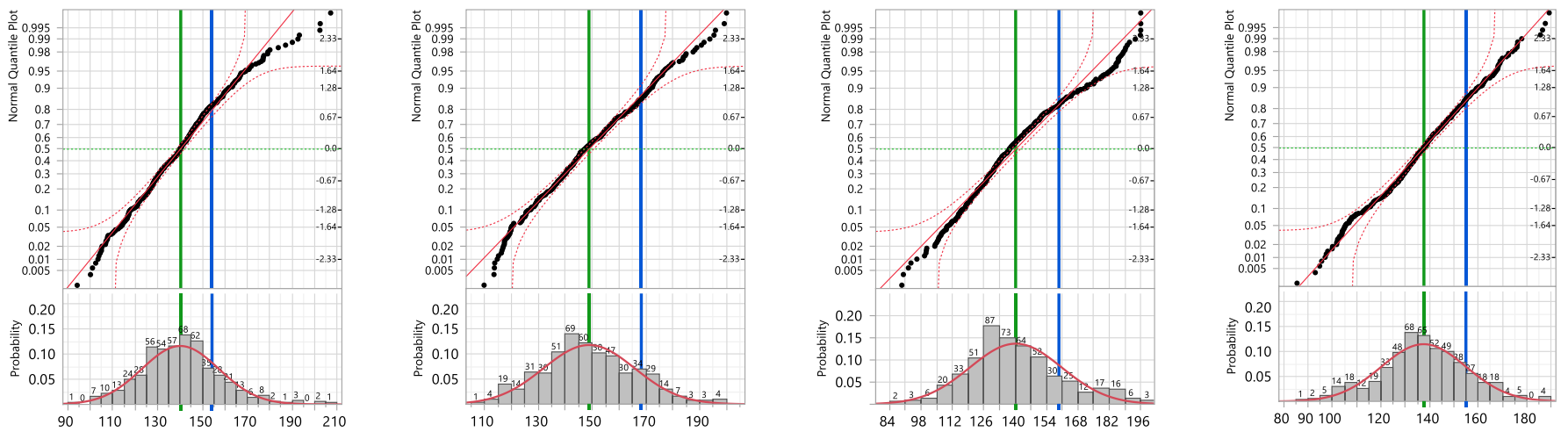

**Design Option 1** **Design Option 2** **Design Option 3** **Design Option 4**

**Fig. 4.** Simulation results in terms of building annual thermal load.

Fig. 4. shows the distribution of results with histograms (bottom), and normality plots (top). In each top graph shown in Fig. 4. , the black points representing the data points are distributed around a diagonal red line. The more the datapoints are falling along the diagonal line, the closer is the distribution to the normal distribution. The horizontal green dashed lines depict the median point in the data set for each design option and the red dashed lines show the confidence limits. The confidence interval is set to 0.95 for all the data sets. The histograms (bottom graphs) show the relative frequency of the results. The vertical green lines show the mean values, and the vertical blue lines represent the deterministic results associated with each design option.

Looking at the plots, there is a clear indication of lack of fit to normal distribution in most design options. Especially, the main difference from normality is evident in the tails rather than in the middle. Furthermore, the normality of the distributions is assessed using the Shapiro Wilk W test (goodness of fit test). In this test, the null hypothesis (H0) is that the data are forming a normal distribution. A small p-value rejects the null hypothesis, meaning there is enough evidence that the data are drawn from a non-normal population. The test results are listed in Table 4.

**Table 4.** Shapiro Wilk W test results.

| | Design Option 1 | Design Option 2 | Design Option 3 | Design Option 4 |
|---|---|---|---|---|
| **W** | 0.985816 | 0.991713 | 0.976497 | 0.995059 |
| **Prob<W** | <.0001 | 0.0069 | <.0001 | 0.1110 |

Note: H0 = The data is from the Normal distribution. Small p-values reject H0.

The null hypothesis for this test is that the data are normally distributed. The Prob < W value listed in the output is the p-value. If the chosen alpha level is 0.05 and the p-value is less than 0.05, then the null hypothesis that the data are normally distributed is rejected. If the p-value is greater than 0.05, then the null hypothesis is not rejected [29]. The results show that the p-values in the first three design options are less than the predefined significance level (0.05). Thus, we can reject the null hypothesis and conclude that the data are not from populations with normal distributions in those design options. The reason could be drawing some input variables from Poisson distribution, also the nonlinear nature of equations in the building energy simulations. On the other hand, in Design Option 4 we cannot reject the null hypothesis, since the p-value is larger than 0.05.

Fig. 5. shows the boxplots to further discuss the probability distributions of the results for each design option. The data points, quantiles, mean values, standard deviations, and deterministic result for each design option are superimposed on the quintile box plot.

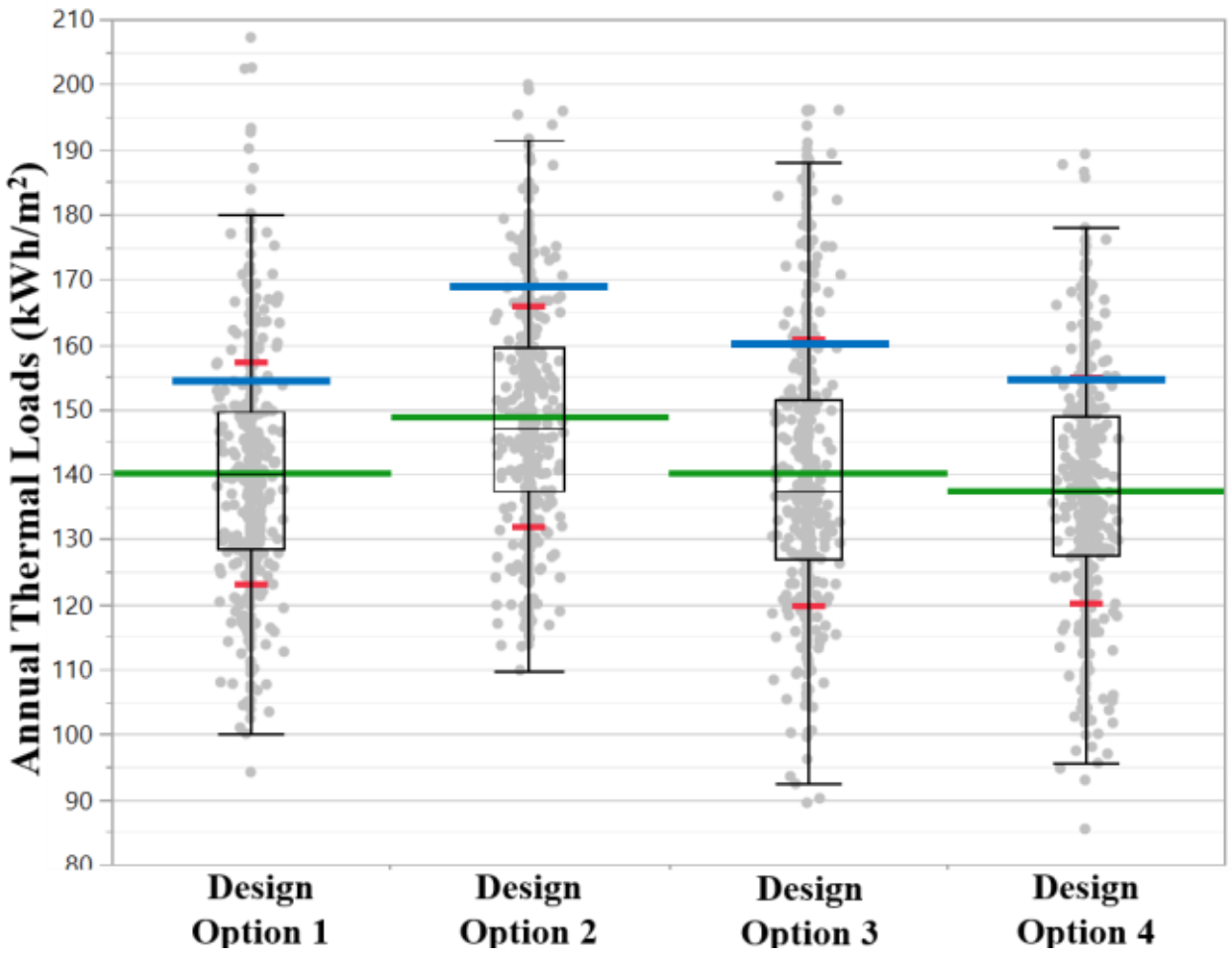

**Fig. 5.** Box Plot of building annual thermal load vs. design type.

The gray points illustrate the data points and the boxplots (shown in black lines) depict the quantiles, dividing the range of the data into four continuous intervals with equal probabilities (25%). The red lines on each boxplot show the standard deviation of the results. The green lines show the mean value for each design option, compared to the blue lines that show the deterministic results. The summary of the results for each design option is described as follows:

1. Design Option 1 shows a range of expected annual thermal load from 94.19 kWh/m$^2$ to 207.17 kWh/m$^2$ with a mean value of 140.16 kWh/m$^2$, a standard deviation of 17.19 kWh/m$^2$, and a variance of 295.50 (kWh/m$^2$)$^2$. The deterministic result shown by the KPI value predicts the building annual thermal load to be equal to 154.42 kWh/m$^2$, which is in the last quartile (located in the fourth 25% of the data).
2. Design Option 2 shows a range of thermal load from 109.85 kWh/m$^2$ to 199.93 kWh/m$^2$ with a mean value of 148.88 kWh/m$^2$, a standard deviation of 16.96 kWh/m$^2$, and a variance of 287.64 (kWh/m$^2$)$^2$. The deterministic result shown by the KPI value is equal to 168.19 kWh/m$^2$, which is in the last quartile (located in the fourth 25% of the data).
3. Design Option 3 shows a range of thermal load from 89.48 kWh/m$^2$ to 196 kWh/m$^2$, with a mean value of 140.29 kWh/m$^2$, a standard deviation of 20.52 kWh/m$^2$, and a variance of 421.07 (kWh/m$^2$)$^2$. The deterministic result indicated by the KPI value shows the value of 159.71 kWh/m$^2$, which is in the last quartile (located in the fourth 25% of the data).
4. Design Option 4 shows a range of thermal load from 85.45 kWh/m$^2$ to 189.22 kWh/m$^2$ with a mean value of 137.58 kWh/m$^2$, a standard deviation of 17.4 kWh/m$^2$, and a variance of 302.76 (kWh/m$^2$)$^2$. The deterministic result shown by the KPI value is equal to 154.76 kWh/m$^2$, which is in the last quartile (located in the fourth 25% of the data).

Building annual thermal load is identified with quantiles in Fig. 6. Using this data, the building thermal loads could be compared under different decision-making criteria.

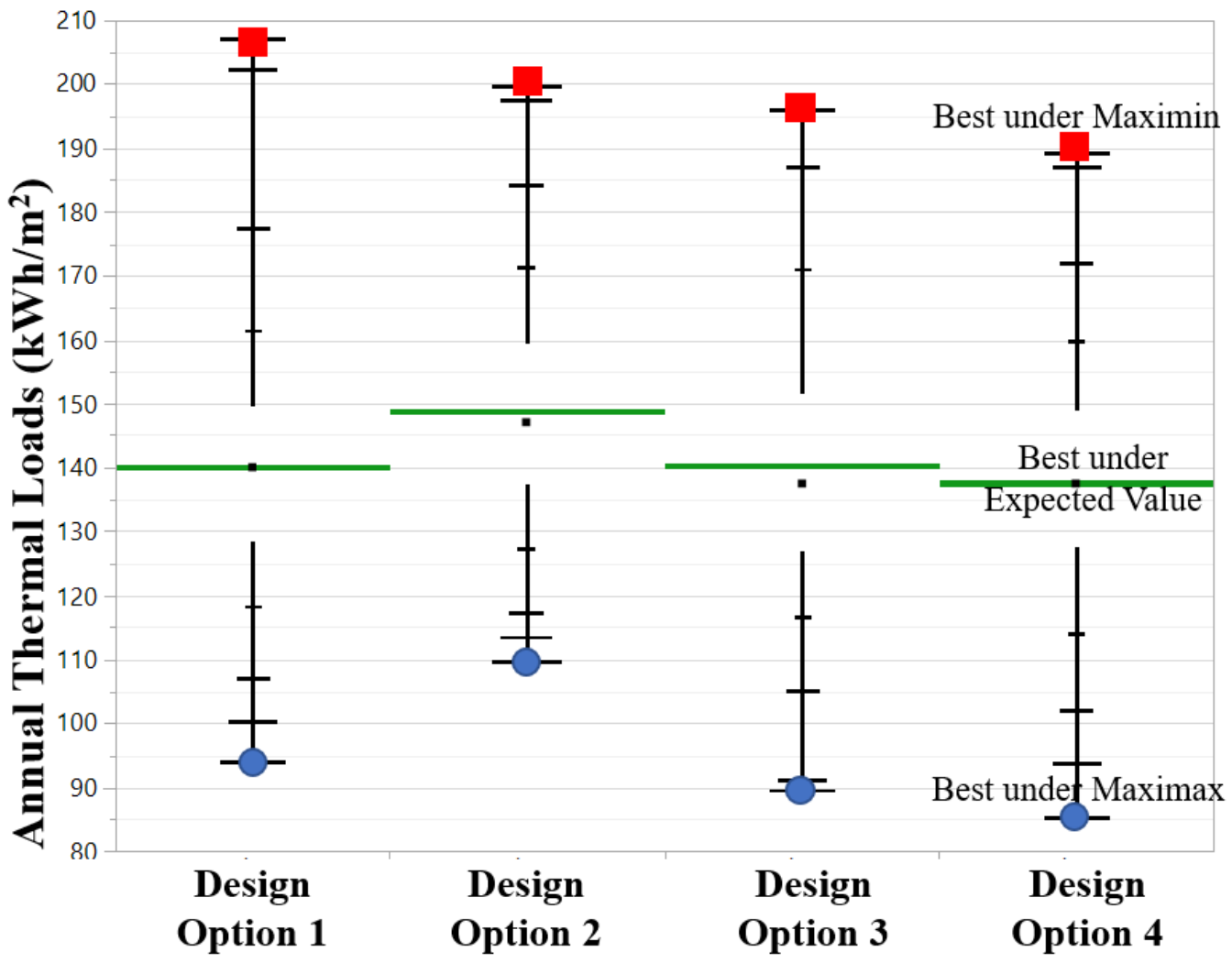

**Fig. 6.** Design decision making suggestions under expected value, maximax, and maximin criteria for test case 3, phase 1.

According to Fig. 6, the minimum and maximum thermal loads are lower in the case of Design Option 4 (85.45 kWh/m$^2$ and 189.22 kWh/m$^2$, respectively) compared to the other design options.

The effect of deterministic and probabilistic results on the ranking of the design options based on different decision-making criteria are summarized as follows:

1. Deterministic: based on the results shown as blue lines in Fig 4., Design Option 1 has the best performance, followed by Design Options 4, 3, and 2.
2. Probabilistic:

   2.1 Expected value criterion: Design Option 4 has the best performance, followed by Design Options 1, 3, and 2.
   2.2 Maximax criterion: Design Option 4 has the best performance, followed by Design Options 3, 1, and 2.
   2.3 Maximin criterion: Design Option 4 has the best performance, followed by Design Options 3, 2, and 1.

Based on the deterministic results (predicted KPIs), it can be concluded that Design Option 1 has the best energy performance, followed by Design Options 4, 3, and 2. However, the probabilistic results provide a more comprehensive outcome with KRIs. Design Option 4 has the lowest mean value of thermal load and is the best design option based on the expected value decision making criterion. Based on the maximax and maximin criteria Design Option 4 is confirmed to have the best performance, while the ranking of other design options is different. The major finding here is that including the uncertainties of input variables in the simulations can lead to probability distributions of the output and may change the performance ranking of design options according to the risk attitudes of design decision makers.

## 6 Conclusion

This research proposes a new framework to implement probabilistic methods in the field of building thermal energy analysis. The proposed framework introduces BIM-ProbE, a new tool (Revit add-in) to create probability distributions for thermal properties of building materials. This framework integrates a building design process with Monte Carlo uncertainty analysis using parametric tools (Grasshopper and its add-ins). This study demonstrated the application of the framework with a test case, considering several sources of uncertainties in building energy analysis with two different probability distributions: normal and Poisson as examples. The test case showed that it takes only about two hours to prepare the data and run the additional simulations to provide probabilistic results and predictions about the possible range of building thermal energy consumption.

The result of this study shows that probabilistic building performance analysis versus deterministic models provides important information about the consequences of design decisions under different possible conditions and may change the performance ranking order of design options. The major finding of this research is that compared with the existing deterministic method for architectural design, using probabilistic methods can result in significantly different design decisions to be made. In addition, using different decision-making criteria including expected value, maximax, and maximin may suggest different design options to select based on different attitudes towards risk. This framework can be widely applied to other design problems and domains to enhance the process of design decision-making. For the future research, the authors intend to work on maximax and maximin criteria based on, for example, first and 99$^{th}$ percentiles instead of the most extreme values in the output range. In this way, inference can be applied to find statistically significant differences between design options.

## Acknowledgement

The authors would like to thank Perkins+Will, who provided the gift funding for this research.